\documentclass[11pt]{amsart}

\newtheorem{theorem}{Theorem}[section]

\theoremstyle{definition}
\newtheorem{definition}[theorem]{Definition}

\theoremstyle{remark}
\newtheorem{remark}[theorem]{Remark}

\numberwithin{equation}{section}

\newcommand{\blankbox}[2]

\begin{document}

\title{New lattice sphere packings denser than Mordell-Weil lattices}

\author{Hao Chen}
\address{Software Engineering Institute,
East China Normal University, Shanghai 200062, P.R. China}
\email{haochen@sei.ecnu.edu.cn}
\thanks{The work was supported by
National Science Foundation of China Grants 10871068 and
11061130539.}



\date{October 28, 2011}


\begin{abstract}1) We present new lattice sphere packings in Euclid spaces of
many dimensions in the range $3332-4096$,  which are denser than known densest Mordell-Weil
lattice sphere packings in these dimensions. Moreover it is proved that if there were some nice linear
binary codes we could construct lattices denser than Mordell-Weil lattices of the many dimensions in the range $128-3272$. \\
2) Lattice sphere packings of many dimensions in the range $4098-8232$ better than present records are presented. Some new dense lattice sphere packings of moderate dimensions $84, 85, 86, 181-189$ denser than any previously known sphere packings of these dimensions are also given.\\
3) New lattices with densities at least 8 times of the densities of
Craig lattices of the dimensions $p-1$, where $p$ is a prime
satisfying $p-1 \geq 1222$, are constructed. Some of these lattices provide new
record sphere packings. \\

The construction is based on the analogues of Craig lattices.

\end{abstract}

\maketitle

\section{Introduction}

The problem to find the dense packing of infinite equal
non-overlapping spheres in Euclid space ${\bf R}^n$ is a classical
mathematical problem(\cite{Kepler,Gauss,CS1,S}). Low dimension sphere packing
problems seem to be understood better than the problems in higher
dimensions. The root lattices in Euclid spaces of dimensions 1, 2,
3, 4, 5, 6, 7 and 8 had been proved  to be the unique densest
lattice sphere packings in these dimensions(see \cite{CS1}). Kepler
conjecture about the densest 3 dimensional sphere packing problem
was proved in \cite{Ha}. Many known densest sphere packings are
lattice packings or packings from finitely many translates of
lattices(see \cite{CS1,CS2, LeS,Va1,Va2}). Constructing lattice
sphere packings from error-correcting codes, algebraic number fields
and algebraic geometry have been proposed by many authors and
stimulated many further
works(\cite{Le,LeS,CS1,Cr1,Cr2,BacNe,Que1,Que2,Que3,Nebe,Ne,Gro,GE, El1,El2,El3,Carm,Flo,Shi1,Shi,Oes,Va1,Va2,Shim1,Shim2}).
Recently Leech lattice, which was found in 1965 in \cite{Le}, has
been proved to
be the unique densest lattice packing in dimension 24 (see \cite{CoEl,CoKu}). For Rogers upper bound and Kabatiansky-Levenshtein upper bound on the
densities of sphere packings we refer to \cite{CS1} pages 19-20. The
recent work \cite{CoEl} gave better upper bounds on densities of
sphere packings. Table 3 in page 711 of \cite{CoEl}
is the latest upper bound for center densities of  sphere packings  in dimensions 1-36.  From Voronoi theory (\cite{Sch,Mar}), there are algorithms to determine the densest lattice sphere packings in each dimension. However the computational task for dimensions $n \geq 9$ is generally infeasible.\\\\

Laminated lattices were known from 1877 (\cite{KZ}) and we refer to \cite{CS1,Conway94}
for a historic survey. In 1982 J. H. Conway and N. J. A. Sloane
published \cite{CS2} in which all densities of laminated lattices up
to dimension 48 were determined. These laminated lattices are known
densest lattices in dimensions 1-29 except dimensions 10, 11, 12, 13
at that time. For dimensions bigger than or equal to $25$, there are many laminated lattices
with the same density. The known densest sphere packing in dimension 12 is
the Coxeter-Todd lattice packing and the known densest sphere
packings in dimensions 10,11,13 were found by Leech and Sloane in
1970 by using non-linear binary codes(see \cite{LeS}). This
situation was changed in 1995 and 1997, when A. Vardy published [22]
and R. Bacher published \cite{Ba}. New non-lattice sphere packings
better than laminated lattice packings were constructed in Euclid
spaces of dimensions 20(see \cite{Va1}), 22(see \cite{CS3}),
27,28,29, 30(see \cite{Va2}) and 18(see \cite{BiEd}). Lattice sphere
packings in Euclid spaces of dimensions 27, 28 and 29 which are
better than laminated lattice packings were constructed in
\cite{Ba}. We refer to Nebe-Sloane list \cite{Nelist} for known densest sphere packings in low dimensions.\\

The knowledge about high dimension sphere packing problem is quite
different. For the high dimensions $n$, in the range $80 \leq n \leq
4096$, $n=2p-2$ where $p$ is a prime number satisfying $p \equiv 5$
$mod$ $6$, or $n=2^t$, where $7 \leq t \leq 12$, the known densest
sphere packings are lattices from algebraic curves over function
fields. That are Mordell-Weil lattices which were discovered  by N. Elkies and
T. Shioda in 1990's (see \cite{El1,El2,Shi1,Shi,Oes}). For Mordell-Weil
lattice in the dimension $n=2p-2$, where $p$ is a prime number
satisfying $p \equiv 5$ $mod$ $6$, the center density is
$\frac{((p+1)/12)^{p-1}}{p^{(p-5)/6}}$ (see \cite{Shi1,Shi} for
other cases). For the center densities of Mordell-Weil lattices in
the dimensions $n=2^t$, where $7 \leq t \leq 12$, we refer to
\cite{El1,CS1} (\cite{CS1} Preface to Third Edition, page xviii).
For example, the known densest sphere packing in dimension $4096$ is
the Mordell-Weil lattice with center density $2^{11527}$. After the discovery of Mordell-Weil lattices by Elkies and Shioda in 1990's, there are a
lot of effort for understanding its structure and construction(for example, \cite{Gro,Nebe,GN}). For
dimensions in the range $149 \leq n=p-1 \leq 3001$ (where $p$ is a
prime number, $149 \leq n \leq 3001$ except $p=509,513$ and $521$),
many of the known densest sphere packings in high $n=p-1$ dimension Euclid
spaces are Craig lattice packings and their recent improvements(see
\cite{Cr1,Cr2,Flo,CS1}). For the sphere packing problem in Euclid spaces
of many dimensions in the range $4100 \leq n \leq 12754608$ or $n
\leq 8 \cdot 10^8$, the best known lattices
were given by Bos-Conway-Sloane construction (\cite{BoCS}, \cite{CS1}, page 17, Table 1.3 and Chapter 8 section 10). The higher dimensional sphere
packing problem was remarked in  N. J. A. Sloane's talk in ICM 1998 as follows." ...But we know very little about this range (of dimensions $80-4096$)...."(\cite{S}).\\

In this paper we propose the analogous Craig
lattice, which is a far-reaching extension of
Criag lattice(\cite{Cr1,Cr2,CS1}, section 6, Chapter 8). These lattices are polynomial lattices
and can be constructed for all
dimensions. Our construction
gives lattice sphere packings of dimensions $n$ in the range $n=p-1
\geq 1222$ better than Craig lattices and their refinements in
\cite{Flo}. Thus some of these lattices provides new record sphere packings
. Our construction also establishes a close
relation between nice lattices and linear error-correcting codes. If
there were some "good enough" linear codes over ${\bf F}_2,{\bf
F}_4,{\bf F}_8$ (see \cite{Gra}), then our construction would lead
to new lattices denser than Mordell-Weil lattices in very many
dimensions in the range $128-3272$. In many high dimensions in the range $3332-4096$, we
prove that some wanted codes in the above description exist. Thus we
present some lattice sphere packings in these dimensions, which are
denser than the Mordell-Weil lattices. New lattices denser than
Shimada lattices in dimensions $84,85$ and $86$(\cite{Nelist,
Shim1,Shim2}) are constructed in section 4.
 Dense lattice sphere packings in many high dimensions in the range $4098-8323$ which are denser
  than the known best lattices from Bos-Conway-Sloane construction are also presented in section 7. Other new dense lattice sphere packings of some moderate dimensions in the range $148-200$ are constructed in section 8.\\

All lattices constructed in this paper are integral lattices.\\

\begin{definition}
For a packing of infinite equal non-overlapping spheres in ${\bf R}^n$ with
centers ${\bf x}_1,{\bf x}_2, ...,{\bf x}_m,....$, the packing radius $\rho$
is defined as $\frac{1}{2}min_{i \neq j}\{||{\bf x}_i-{\bf
x}_j||\}$. The density $\Delta$ is  $lim_{t \rightarrow 0}
\frac{Vol\{{\bf x} \in {\bf R}^n: ||{\bf x}|| <t, \exists {\bf x}_i,
||{\bf x}-{\bf x}_i||<\rho\}}{Vol\{{\bf x} \in {\bf R}^n: ||{\bf
x}||<t\}}$. Center density $\delta$ is defined as
$\frac{\Delta}{V_n}$ where $V_n$ is the volume of the ball of
radius $1$ in ${\bf R^n}$.
\end{definition}

Let ${\bf b}_1,...,{\bf b}_m$ be $m$ linearly independent vectors in
the Euclid space ${\bf R}^n$ of dimension $n$. The discrete point
sets ${\bf L}=\{x_1{\bf b}_1+\cdots+x_n{\bf b}_m: x_1,...,x_m \in
{\bf Z}\}$ is rank $m$ lattice in ${\bf R}^n$.  The determinant of
the lattice is defined as $det({\bf L})=det(<{\bf b}_i, {\bf b}_j>$.
The volume of the lattice is $Vol({\bf L})=(det({\bf
L})^{\frac{1}{2}}$. Let $\lambda({\bf L})$ be the Euclid norm of the
shortest non-zero vectors in the lattice and the minimum norm of the
lattice is just $\mu({\bf L})=(\lambda({\bf L}))^2$. When the
centers of the spheres are these lattice vectors in ${\bf L}$ we
have $\rho=\frac{1}{2}\lambda({\bf L})$ and the center density
$\delta({\bf L})=\frac{\rho^n}{Vol({\bf L})}$. The lattice ${\bf
L}^{*}=\{{\bf y} \in {\bf R}^m: <{\bf y},{\bf x}> \in {\bf Z}\}$ is
called the dual lattice of the lattice ${\bf L}$. A lattice is
called integral if the inner products between lattice vectors are
integers. An unimodular lattice is the integral lattice ${\bf L}$
satisfying ${\bf L}^{*}={\bf L}$. The unimodular lattices whose
minimum norm attain the largest possible bound are called extremal
unimodular lattices. The problem to find extremal unimodular
lattices seems very difficult and  has attracted many works(see
\cite{BacNe,Ne,CS1}, Cha.7 or Nebe-Sloane database of
lattices \cite{Nelist}).\\

Let $r$ be a prime power and ${\bf F}_r$ be the finite field with
$r$ elements. A linear (non-linear ) error-correcting code ${\bf C}
\subset {\bf F}_r^n$ is a $k$ dimensional subspace(or a subset of
$M$ vectors). For a codeword ${\bf x} \in {\bf C}$, $wt({\bf x})$ is
the number of nonzero coordinates of ${\bf x}$. The minimum Hamming
weight(or distance) of the linear (or non-linear) code $C$ is
defined as $d({\bf C})=min_{{\bf x}\neq {\bf y},{\bf x},{\bf y} \in {\bf C} }\{wt({\bf x-y})\}$.
We refer to $[n,k,d]$ (or $(n,M,d)$) code as linear (or
non-linear)code with length $n$,  distance $d$ and dimension $k$ (or
$M$ codewords). Given a binary code ${\bf C} \subset {\bf F}_2^n$
the construction A (\cite{CS1,LeS}) leads to a lattice in ${\bf
R}^n$. The lattice ${\bf L(C)}$ is defined as the set of integral vectors ${\bf
x}=(x_1,...,x_n) \in {\bf Z}^n$ satisfying $x_i\equiv c_i$ $mod$ $2$
for some codeword ${\bf c}=(c_1,...,c_n) \in {\bf C}$. It is east to
check $\rho=\frac{1}{2}min\{\sqrt{ d({\bf C})},2\}$ and $Vol({\bf
L(C)})=2^{n-k({\bf C})}$ where $k({\bf C})$ is the dimension of the
code ${\bf C}$. This gives a lattice sphere packing with the center
density $\delta=\frac{min\{\sqrt{ d({\bf C})},2\}^n}{2^{2n-k({\bf
C})}}$. This construction A leads to some best known densest lattice
packings in low dimensions(see \cite{CS1}). For a non-linear binary
code, the same construction gives the non-lattice packing with
center density $\frac{M \cdot min\{\sqrt{d({\bf
C})},2\}^n}{2^{2n}}$, where $M$ is number of codewords in the
nonlinear binary code ${\bf C}$. For example, the non-linear length
$11$, minimum distance $4$ binary codes with $72$ codewords implies
a non-lattice packing with center density
$\delta=\frac{9}{256}=0.03516$, which is the known densest sphere
packing in dimension 11(\cite{LeS,CS1}).  The known densest sphere
packings in Euclid spaces of dimensions 10 and 13 were constructed
similarly in \cite{LeS}. The packing in dimension 10 is from
non-linear binary $(10, 40, 4)$ code(\cite{CS1}). From the
non-linear code $(12,144,4)$ a non-lattice sphere packing with
center density $\frac{9}{256}$ (which is smaller than the center
density $\frac{1}{27}$ of the Coxeter-Todd lattice) can be
constructed (\cite{LeS,CS1}). The packing in dimension 13 is the
union of infinitely many copies of this non-lattice dimension 12
packing from non-linear $(12, 144, 4)$ binary codes(see
\cite{LeS}). We refer to \cite{Nelist} for records of dense sphere packings in Euclid spaces of various dimensions.\\

\section{Analogous Craig lattices}

Let $\zeta$ be a primitive $p$-th root of unity, where $p$ is an odd
prime. Then the ring of the algebraic integers in ${\bf Q}[\zeta]$
is ${\bf Z}[\zeta]$(\cite{Cr1,Cr2}). The Craig lattice ${\bf
A}_{p-1}^{(i)}$ of rank $p-1$ introduced in \cite{Cr1} is the ideal
in the ring ${\bf Z}[\zeta]$(free ${\bf Z}$ module) generated by
$(1-\zeta)^i$, where $i$ is a positive integer. A cyclotomic construction of Leech lattice was given by  Craig in \cite{Cr2}. In \cite{Flo}
refinements of Criag lattices  were constructed by adding some
fractional numbers. The refinements of Craig lattices have their
center density at least the three times of the center density of
the original Craig lattices. These provided some new record lattice sphere packings in some dimensions(see \cite{Flo}).\\

Another form of Craig lattices was given in \cite{CS1} section 6 of
Chapter 8. It is  a cyclic lattice ${\bf A}_{n}^{(m)}$ of rank $n$ in the
ring ${\bf Z}[x]/(x^{n+1}-1)$, the ideal generated by
$(x-1)^m$ in the ring ${\bf Z}[x]/(x^{n+1}-1)$. The volume of the
Craig lattice ${\bf A}_n^{(m)}$ is $(n+1)^{m-\frac{1}{2}}$. When
$n+1=p$ is an odd prime and $m < \frac {p}{2}$, the minimum norm
$\mu({\bf A}_n^{(m)}) \geq 2m$ (see Theorem 7 at page 223 of
\cite{CS1}).
This is just the original Craig lattices introduced in \cite{Cr1}.\\

We propose the following analogue of the Craig lattice as a ${\bf
Z}$ sub-module of the ${\bf Z}$ module $R={\bf Z}<1,x,...,x^n>$
spanned by $1, x, x^2,...,x^n$. Given positive integers $n, m,l$
satisfying $m < \frac{n}{2}$ and $ l \geq n+1$, consider the
sub-module ${\bf A}_n^{(m,l)}={\bf Z} (x-1)^n+{\bf
Z}(x-1)^{n-1}+\cdots +{\bf Z} (x-1)^m + l{\bf Z} (x-1)^{m-1}+l{\bf
Z}(x-1)^{m-2}+\cdots+l{\bf Z}(x-1)$ in ${\bf Z}<1,x,...,x^n>$. Then
any element $v$ in ${\bf A}_n^{(m,l)}$ can be represented as $v_0
1+v_1 x+ \cdots+v_n x^n$. The set of all coordinates
$(v_0,v_1,...,v_n)$ of these vectors $v$ in ${\bf A}_n^{(m,l)}$ is a
sub-lattice in ${\bf Z}^{n+1}$. Equivalently we take the base
$1,x,...,x^n$ as the orthogonal base of the ${\bf Z}$ module $R$ and
the ${\bf Z}$ sub-module described as above is just the lattice. For any polynomial $f(x)=a_0+a_1x+\cdots+a_nx^n \in R$ with integral
coefficients satisfying $f(1)=0$, $f(x)$ is a linear combination of
$(x-1),...,(x-1)^n$ with integral coefficients. If $l$ has no prime
factor smaller than $m$, a polynomial $f(x)=a_nx^n+\cdots+a_1x+a_0$
with integral
coefficients is in the lattice ${\bf A}_n^{(m,l)}$ if and only if $f(1)=0$ and $f^{(i)}(1)\equiv 0$ $mod$ $l$ for $i=1,...,m-1$.\\

{\bf Theorem 2.1.} {\em  \begin{enumerate}
\item When $n+1$ is a prime or $m=1$, ${\bf
A}_n^{(m,n+1)}$ is just the Craig lattice;

\item Given positive integers $n, m,l$ satisfying $m < \frac{n}{2}$
and $ l \geq n+1$, the ${\bf A}_n^{(m,l)}$ is a lattice of rank $n$
and its volume is $l^{m-1} (n+1)^{\frac{1}{2}}$. When $l$ is a prime number, the minimum norm of the lattice ${\bf
A}_n^{(m,l)}$ satisfies $\mu({\bf A}_n^{(m,l)}) \geq 2m$.
\end{enumerate}}

\begin{proof} 1) When $m=1$, ${\bf A}_n^{(1,l)}$ has an integral base
$\{(x-1), (x-1)x,...,(x-1)x^{n-1}\}$. Thus ${\bf A}_n^{(1,l)}={\bf
A}_n^{(1)}$.\\

We prove the ${\bf A}_n^{(m,l)}$ is a cyclic lattice when $l=n+1$ is
a prime number. It is clear that $(x-1)^jx=(x-1)^{j+1}+(x-1)^j$ is
in ${\bf A}_n^{(m,l)}$ for $m \leq j \leq n-1$ and
$l(x-1)^jx=l(x-1)^{j+1}+l(x-1)^j \in {\bf A}_n^{(m,l)}$ for $1 \leq j
\leq m-1$. We only need to check the element $(x-1)^n
x-(x^{n+1}-1)$, which is a shift of $(x-1)^n$, is in ${\bf A}_n^{(m,l)}$. Since
$(x-1)^nx-(x^{n+1}-1)=(x-1)^{n+1}-x^{n+1}+1+(x-1)^n$, the
coefficients of $(x-1)^j$, where $1 \leq j \leq n-1$, in the expansion of $(x-1)^{n+1}-x^{n+1}+1$ can be divided by $n+1$,
when $n+1$ is a prime number. Thus $(x-1)^n x-(x^{n+1}-1)$ is also in ${\bf A}_n^{(m,l)}$ when $n+1$ is
a prime. Then ${\bf A}_n^{(m,n+1)}$ contains the Craig lattice ${\bf
A}_n^{(m)}$ as a sub-lattice. It has the same index as
a sub-lattice of the lattice ${\bf A}_n^{(1)}$. The conclusion is
proved.\\
2) It is obvious that ${\bf A}_n^{(m,l)}$ is a sub-lattice of the
rank $n$ Craig lattice ${\bf A}_n^{(1)}=\{(v_0,v_1,...,v_n) \in {\bf
Z}^{n+1}: v_0+v_1+\cdots+v_n=0\}$. On the other hand ${\bf
A}_n^{(m,l)}$ has index $l^{m-1}$ in the Craig lattice ${\bf
A}_n^{(1)}$. Thus it is a rank $n$ lattice and has volume $l^{m-1}
Vol({\bf A}_n^{(1)})=l^{m-1}(n+1)^{\frac{1}{2}}$.
The proof of the second conclusion is the same as the proof of
Theorem 7 of page 223 of \cite{CS1}. If $\mu({\bf A}_n^{(m,l)}) <2m$ We have an element $f(x)=\Sigma_{i \in S}
x^i-\Sigma_{j \in T} x^j \in {\bf A}_n^{(m,l)}$, where $S$ and $T$
are two sub-sets in $\{0,1,...,n\}$ satisfying $h=|S|=|T| <m$. Here
$S=\{s_1,...,s_h\}$ and $T=\{t_1,...,t_h\}$ may contain repeated
elements. Then from the condition  $f(x) \in {\bf A}_n^{(m,l)}$ we
have $f^{(i)}(1) \equiv 0$ $mod$ $l$, for $i=0,1,...,m-1$. Then we
have $\Sigma_{j=1}^h s_j^i \equiv \Sigma_{j=1}^h t_j^i$ $mod$ $l$,
for $i=0,1,...,m-1$. Since $l$ is a prime number,
from the Newton's identities over the finite field ${\bf Z}/l{\bf Z}$, the elementary symmetric functions of
$S$ and $T$ of degree $<m$ have to be the same. Thus $S=T$ and
$f(x)=0$ since we have $ l \geq n+1$.
\end{proof}
From Theorem 2.1 we have the following analogous Craig lattices for
all  dimensions  with "not bad" densities. The following result can be compared with the construction
of Craig-like lattices in \cite{Carm}. It is obvious our lattices are denser, since in their construction
the prime number $q$ is required to be the smallest prime $q$ satisfying $ q \equiv 1$ $mod$ $n$.\\

{\bf Theorem 2.2.} {\em For each dimension $n$ and each $m <\frac{n}{2}$ we have a analogous
Craig lattice with density $\Delta_n \geq
\frac{m^{\frac{n}{2}}}{2^{m-1+\frac{n}{2}}\cdot
n^{m-1}(n+1)^{\frac{1}{2}}}$. When suitable $m$ is taken we have
$\frac{1}{n}log_2 \Delta_n \geq -\frac{1}{2}log_2log_2 n+o(1)$.}\\

\begin{proof}
From the Bertrand postulate there exists a prime number $l$ between
$n$ and $2n$ for any positive integer $n$. Set $m$ to be the integer
nearest to $\frac{n}{2 log_e n}$ as in \cite{CS1} page 224. We take
the analogue Craig lattice ${\bf A}_n^{(m,l)}$. A direct calculation
gives us the result.
\end{proof}

Let $l$ be an odd number. We define a mapping $\pi: {\bf A}_n^{(m,l)}/2({\bf A}_n^{(m,l)}) \rightarrow {\bf Z}^{n+1}/2({\bf Z}^{n+1})$
by $\pi(a_0+a_1x+\cdots+a_nx^n)=(a_0,...,a_n)$ $mod$ $2$. This  is a injective ${\bf Z}$ linear mapping.
It is clear that $2(\frac{f^{(i)}(1)}{i!})$ can be divided by $l$ implies that $\frac{f^{(i)}(1)}{i!}$ can be divided by $l$, since $l$
is an odd number. We can check
that the image of $\pi$ is the linear binary $[n+1,n,2]$ code.\\

{\bf Theorem 2.3.} {\em Suppose the positive integers $n, m,l$ satisfy $m <
\frac{n}{2}$, $ l \geq n+1$ and $l$ is a odd prime number.  If there exists a linear binary
sub-code of the $[n+1,n,2]$ code with parameters $[n+1, k, \geq
8m]$, then we have a lattice with center density at least
$\frac{2^{k-\frac{n}{2}} \cdot m^{\frac{n}{2}}}{l^{m-1}
(n+1)^{\frac{1}{2}}}$.}\\

\begin{proof}
The $[n+1,k, \geq 8m]$ binary linear sub-code $V$ is
in the image $\pi({\bf A}_n^{(m,l)})$ as a binary linear
$[n+1,n,2]$ code. From the ${\bf Z}$-linearity of $\pi$, the
inverse image $\pi^{-1}(V)$ is a lattice with volume $2^{n-k}
Vol({\bf A}_n^{(m,l)})$. Let ${\bf v}$ be a vector in
$\phi^{-1}(V)$. If $\phi({\bf v})=0$, ${\bf v} \in 2{\bf
A}_n^{(m,l)}$, then the Euclid norm of ${\bf v}$ is at least $8m$.
If $\pi({\bf v}) \neq 0$, then at least $8m$ coordinates of the
vector  ${\bf v}$ are odd numbers and the Euclid norm of ${\bf v}$
is at least $8m$. The conclusion is proved.
\end{proof}

{\bf Theorem 2.4.} {\em Suppose the positive integers $n, m,l$ satisfy $m \leq
\frac{n+1}{2}$, $ l \geq n+1$ and $l$ is a odd prime.  If there exists a linear binary
$[n,k,8m]$ code then we have a lattice with center
density at least $\frac{2^{k-\frac{n}{2}} \cdot
m^{\frac{n}{2}}}{l^{m-1} (n+1)^{\frac{1}{2}}}$.}\\

\begin{proof}
The extended code of the binary $[n,k,8m]$ code by adding a parity
check column we get the linear sub-code in the Theorem 2.3.
\end{proof}

From Theorem 2.4 some lattices of rank $n=p^2-1$,
where $p$ is a prime, can be constructed. We have a lattice ${\bf A}_{120}^{(11,127)}$ of rank $120$ with center density at least
$2^{74.0640}$(less than the center
density $2^{76}$ of the Bos-Conway-Sloane construction $\eta({\bf
E}_8)$, see \cite{CS1} p.242), a lattice ${\bf A}_{168}^{(13,173)}$ of rank $168$ with center
density at least $2^{133.9011}$(larger
than the center density $2^{120}$ of the Bos-Conway-Sloane
construction $\eta({\bf \Lambda}_{24})$, see \cite{CS1} p.242),  a
lattice ${\bf A}_{288}^{(17,293)}$ of rank $288$ with center density at least
$2^{309.3031}$(larger than $2^{300}$
of the Bos-Conway construction $\eta({\bf \Lambda}_{24})$, see
\cite{CS1} p.242), a lattice ${\bf A}_{360}^{19,367)}$  of rank $360$ with center density at
least $2^{427}$(larger than
$2^{408}$ of the $\eta({\bf \Lambda}_{24})$ in ${\bf R}^{360}$, see
\cite{CS1} p.242). \\

We now apply Theorem 2.4 to construct lattices in dimensions  $60$,
$96$, $136$ and $144$. In the dimension 60, the two known good lattices
are Kashichang-Pasupathy lattice ${\bf Ks}_{60}$ with the center
density $2^{17.4346}$(see \cite{CS1}, page xivi) and the extremal $60$
dimensional lattice with center density $(\frac{3}{2})^{30} \approx
2^{17.55}$(see Nebe-Sloane database of lattices \cite{Nelist}). Applying
Theorem 2.4 to ${\bf A}_{60}^{(7,61)}$
and the binary linear $[60,1,60]$ code(see
\cite{Gra}) we get a lattice with center density $2^{16.672}$(also from \cite{Flo}). From
\cite{Gra} if there was a binary linear $[60,27,16]$ code, then a
new denser lattice with center density $2^{18.1039}$ could be
constructed. The Elkies lattice in dimension $60$ has center density $2^{19.04}$, which is the section
of Mordell-Weil lattice(see \cite{Va2}). If there was a binary linear $[60, 28, 16]$ code(\cite{Gra}), we could construct
a lattice sphere packing in dimension $60$ with center density $2^{19.1039}$. The known densest sphere packing in dimension $96$ is  is the lattice
$\eta({\bf P}_{48q})$ with center density $2^{52.078}$(see \cite{CS1} p.16). From the table in
\cite{Gra}, there exists a linear binary $[96,23,32]$ code, we get a lattice with center density
$2^{47.9003}$. If there was a linear binary
$[96, 30, 32]$ code(\cite{Gra}), we could construct a lattice sphere packing
with center density $2^{54.9003}$. The known densest sphere packing in dimension $136$ is the lattice
$\eta({\bf E}_8)$ with center density $2^{100}$(see \cite{CS1} p.16)
from Bos-Conway-Sloane construction(\cite{BoCS}). Applying Theorem
2.4 to ${\bf A}_{136}^{(4,137)}$ and binary linear $[136,47,32]$
code (see \cite{Gra}), we have a lattice with center density
$2^{90.1570}$. However from \cite{Gra} it is possible that there was a
$[136,57,32]$ code, which would lead to a possible new lattice with
center density $2^{100.1570}$. In dimension $144$ there is a dense
lattice $\eta({\bf \Lambda}_{24})$ with center density $2^{96}$(see
\cite{CS1} p.242) from Bos-Conway-Sloane construction \cite{BoCS}.
From the ${\bf A}_{144}^{(14,149)}$ and the binary linear
$[144,1,144]$ code(see \cite{Gra}) we get a new lattice with center
density $2^{105.6736}$.\\

In dimension $n=160$, since $n+1=161=23 \cdot 7$ is not a prime, we have no Craig lattice
in this dimension. The analogous Craig lattice ${\bf A}_{160}^{(16, 163)}$ has its center density $\delta_{160}=\frac{8^{80}}{163^{15.5}} \approx 2^{126.4051}$. By
using the trivial linear binary $[160, 1, 160 ]$ code and Theorem 2.4 we get a dense lattice in dimension $160$ with center density at least
$2^{127.4051}$. The analogous Craig lattice ${\bf A}_{160}^{(8, 163)}$ has center density $2^{104.8847}$. Applying Theorem 2.4 to the
linear binary code $[160, 19, 64]$(\cite{Gra}) then a lattice with center density $2^{123.8847}$ can be constructed. If there was a linear binary $[160, 27, 64]$ code(\cite{Gra}), a lattice with center density $2^{131.8847}$ could be constructed. On the other hand, there is no Mordell-Weil lattice in dimension 160. The nearest Mordell-Weil lattice in dimension smaller than $160$ (from Theorem 1.1 of \cite{Shi}) is the Mordell-Weil  lattice of dimension 140 with center density $2^{113.31}$. There is no child lattice $\eta({\bf E}_8)$ in dimension $160$(\cite {CS1}, page 241). The lattice in dimension $160$ from Minkowski-Hlawka Theorem has center density approximately $111.2378$. Thus our construction from analogous Craig lattice gives a new dense lattice in dimension $160$.
\remark
When $l(>n)$ is a prime number, the analogous Craig lattice ${\bf A}_n^{(m,l)}$ is just the section of the Craig lattice ${\bf A}_{l-1}^{(m,l)}$ by imposing the condition that the last $l-n-1$ coordinates are zero.
\endremark

\section{Improving Craig lattices and their refinements}

The main result of this section is the following theorem.\\

{\bf Theorem 3.1.} {\em Let $p$ be a prime larger than or equal to $1223$. Suppose ${\bf
A}_{p-1}^{(m)}$ is the densest Craig lattice of dimension $p-1$. We
can construct a new lattice with center density at least
$8\delta({\bf A}_{p-1}^{(m)})$ from Theorem 2.4.}\\

\begin{proof}
It is known the Craig lattice ${\bf A}_{n}^{(m)}$, where $m$ is
nearest integer of $\frac{n}{2 log_e (n+1)}$, is the densest Craig
lattice in the dimension $n=p-1$, where $p$ is a prime number. Since
$n \geq 1222$, then $\frac{8(\frac{n}{2 log_e (n+1)}+1)}{n} \leq
\frac{4}{7}$. Concatenating linear
$[\lfloor\frac{n}{7}\rfloor,1,\lfloor\frac{n}{7}\rfloor]$ code over
${\bf F}_8$ with binary linear $[7,3,4]$ code and a suitable trivial
extension we get  a binary linear $[n, 3 , 8(\frac{n}{2 log_e
(n+1)}+1)]$ code. From Theorem 2.4 we get the conclusion.
\end{proof}

In \cite{Flo}, the Craig lattices are refined to
new lattices with center density at most $6\delta({\bf A}_n^{(m)})$
in the range $1298 \leq n \leq 3482$. Thus our constructed lattices
are better than the lattices in \cite{Flo}. Some of these new dense lattice sphere packings are better than any previously known ones.\\

\begin{table}[ht]
\caption{}\label{eqtable}
\renewcommand\arraystretch{1.5}
\noindent\[
\begin{array}{|c|c|c|c|}
\hline
dim=n&new-log_2 \delta&known-densest&nearest-MW(dim<n)\\
\hline
1398&2908.8254&2905.8254(Craig)&2919.8743(1364)\\
\hline
1432&2980.6910&2977.6910(Craig)&3012.6846(1400)\\
\hline
2178&5131.4554&5128.4554(Craig)&5086.8746(2120)\\
\hline
2296&5592.5709&5589.5709(Craig)&5377.7840(2216)\\
\hline
\end{array}
\]
\end{table}

\section{Lattices in dimensions 52, 68, 84, 85, 86, 120, 168, 242, 246, 248, 288 and 360 from analogous Craig lattices}

The known densest lattices in dimension $48$ and $56$ are extremal
unimodular lattices(\cite{CS1} Preface to Third Edition) and the
known extremal unimodular lattices in dimension $80$ has center
density $2^{40}$ which is slightly lessen than the center density
$2^{40.14}$ of the known densest lattice(Mordell-Weil lattice) in
this dimension(see \cite{BacNe}). The recently constructed extremal
unimodular lattice (see \cite{Ne}) in dimension $72$ is the known
densest lattice with center density $2^{36}$. Therefore we compare
the new lattices from Theorem 2.4 with
known (extremal) unimodular lattices. The Gaborit extremal unimodular lattice in dimension $52$ has its minimum norm $5$ and center density $(\frac{5}{4})^{26} \approx 2^{8.4552}$(\cite{Nelist}). Applying Theorem 2.4 to the analogous Craig lattice ${\bf A}_{52}^{(6,53)}$ and the $[52,1,52]$ linear binary code we can construct a integral lattice with center density $2^{10.7045}$. This is a new dense sphere packing in dimension $52$. The
dimension $68$ extremal unimodular lattices with minimum norm $6$ have
their center densities $(\frac{3}{2})^{34} \approx 2^{19.89}$.
Applying Theorem 2.4 to analogous Craig lattice ${\bf
A}_{68}^{(4,71)}$ and binary linear $[68,8,32]$ code(see \cite{Gra})
we get a new lattice with center density at least $2^{20.4757}$. The
volume of this new lattice is $2^{60}\cdot 71^3
\cdot 69^{\frac{1}{2}}$. The Thompson-Smith unimodular lattice in dimension 248 has minimum
norm $10$ or $12$. Thus its center density is at most
$3^{124}\approx 2^{196.54}$(see Nebe-Sloane database of lattices,
\cite{Ne}). Applying Theorem 2.4 to ${\bf A}_{248}^{(4,251)}$ and
binary linear $[248,131,32]$ code (see \cite{Gra}) we get a new
lattice with center density at least $2^{227.0997}$. The volume of
this new lattice is $2^{137}\cdot 251^3 \cdot 249^{\frac{1}{2}}$. In dimension $240$ the known densest sphere packing is the Craig lattice packing with center density $2^{245.0006}$. Applying Theorem 2.4 to the analogous Craig lattice ${\bf A}_{242}^{(12, 251)}$ with center density at least $2^{221.1127}$ and the $[242,22,96]$ linear binary code, we get a dense lattice sphere packing in dimension $242$ with center density $2^{243.1127}$. Applying Theorem 2.4 to the analogous Craig lattice ${\bf A}_{246}^{(12,251)}$ with center density $2^{226.2827}$ and linear binary $[246,23,96]$ code(\cite{Gra}) we get a dense lattice sphere packing in dimension $246$ with center density $2^{249.2827}$.\\

In dimension 104, the presently known densest sphere packing is the Mordell-Weil lattice with center density $\frac{(\frac{9}{2})^{52}}{53^8} \approx 2^{67.0168}$(\cite{Shi}, Theorem 1.3). Before the invention of Mordell-Weil lattices by Elkies and Shioda, the previously known densest sphere packing is the child lattice $\eta({\bf E}_8)$ with center density $2^{60}$(\cite{CS1}, page 17, Table 1.3). The analogous Craig lattice ${\bf A}_{104}^{(5, 107)}$ has center density $2^{38.3949}$. Thus if there was a linear binary $[104, 22, 40]$ code{\cite{Gra}), a lattice denser than $\eta({\bf E}_8)$ with center density
$2^{60.3949}$ could be constructed. However even the codes attaining the upper bound in the Table \cite{Gra} exist, we cannot construct lattices denser than Mordell-Weil lattice in the dimension $104$.\\

In \cite{Shim1,Shim2} long computation of algebraic geometry over finite fields was used to constructed
dense lattice sphere packings in dimensions 84,85 and 86 with center densities $\delta_{84}^{Shimada} \approx 2^{30.795}, \delta_{85}^{Shimada} \approx 2^{32.5}$
and $\delta_{86}^{Shimada} \approx 2^{34.2075}$(see Nebe-Sloane list \cite{Nelist}, and the comment about Shimada's 86 dimensional lattice there). We take $n=84,85,86$ and $l=89$, which is a prime, and $m=4$.
The corresponding dimension $n$ analogous Craig lattice has center density $\delta_n \approx 2^{n/2-22.67}$. Since $[84,16,32],[85,16,32]$
and $[87,17,32]$ linear binary codes exist(\cite{Gra}), we have better lattices in dimension 84 with center density at least $2^{35.4}$,
better lattices in dimension $85$ with center density at least $2^{35.83}$, and better lattice in dimension $86$ with center density at least $2^{37.33}$,
from Theorem 2.4. The analogous Craig lattice ${\bf A}_{86}^{(10,89)}$ has center density $2^{38.3225}$, the analogous Craig lattice ${\bf A}_{85}^{(10,89)}$
has center density $2^{37.1616}$ and the analogous Craig lattice ${\bf A}_{84}^{(10,89)}$ has center density $2^{36.006}$. Applying Theorem 2.4 to the
 trivial linear binary $[84,1,84],[85,1,85]$ and
$[86,1,86]$ codes we get new lattices in dimensions $84, 85$ and $86$ with center densities $2^{37.006}$, $2^{38.1616}$ and $2^{39.3225}$.\\

Applying Theorem 2.4 to the analogous Craig lattice ${\bf A}_{120}^{(11,127)}$ and linear binary $[120,1,120]$ code we get a lattice sphere packing in dimension $120$ with center density $2^{75.0640}$. Applying Theorem 2.4 to the analogous Craig lattice ${\bf A}_{168}^{(13,173)}$ and linear binary $[168,2,114]$ code (\cite{Gra}) we get a lattice sphere packing in dimension $168$ with center density a lattice of rank $168$ $2^{135.9011}$.  From \cite{Gra}, there is a linear binary $[144,9,68]$ code. Thus we have a linear binary $[288, 9,136]$ code. Applying Theorem 2.4 to the analogous Craig lattice ${\bf A}_{288}^{(17,293)}$ and linear binary $[288,9,136]$ code  we get a lattice sphere packing in dimension $288$ with center density $2^{318.3031}$. From \cite{Gra}, there is a linear binary $[180,16,78]$ code. Thus we have a linear binary $[360, 16,156]$ code. Applying Theorem 2.4 to the analogous Craig lattice ${\bf A}_{360}^{(19,367)}$ and linear binary $[360,16,156]$ code  we get a lattice sphere packing in dimension $360$ with center density $2^{443}$.\\

In table 2 we list some new lattice sphere packings from Theorem
2.4.
\begin{table}[ht]
\caption{}\label{eqtable}
\renewcommand\arraystretch{1.5}
\noindent\[
\begin{array}{|c|c|c|}
\hline
dimension&new-log_2 \delta&known\\
\hline
52&10.7045&8.4552(Gaborit), 10.4578(MW)\\
\hline
60&16.672&19.04(Elkies,section-MW), 17.55(extremal)\\
\hline
68&20.6757&19.89(Gaborit+Harada-Kitazume)\\
\hline
84&37.006&30.795(Shimada)\\
\hline
85&38.1616&32.5(Shimada)\\
\hline
86&39.3225&34.2075(Shimada)\\
\hline
96&47.9003&52.078(\eta({\bf P}_{48q}))\\
\hline
120&75.0640&76(\eta({\bf E}_8))\\
\hline
144&105.6736&96(\eta({\bf \Lambda}_{24}))\\
\hline
160&127.4051&111(Minkowski-Hlawka)\\
\hline
168&135.9011&120(\eta({\bf \Lambda}_{24}))\\
\hline
246&249.2827&234.33039(Minkowski-Hlawka)\\
\hline
248&227.0997&196.54(Thompson-Smith)\\
\hline
288&318.3031&300(\eta({\bf \Lambda}_{24}))\\
\hline
360&443&408(\eta({\bf \Lambda}_{24}))\\
\hline
\end{array}
\]
\end{table}

In the following table 3 we list some possible better lattices under
the condition that some nice codes exist. The Elkies lattices of dimensions $57-60$ are cross-sections of Mordell-Weil lattices, we refer to \cite{Va2}, page 278.\\

\begin{table}[ht]
\caption{}\label{eqtable}
\renewcommand\arraystretch{1.5}
\noindent\[
\begin{array}{|c|c|c|c|}
\hline
dimension&possible-log_2 \delta&known-densest&condition\\
\hline
57&16.1040&15.37(Elkies,section-MW)&[57,25,16]\\
\hline
58&17.1040&16.46(Elkies,section-MW)&[58,26,16]\\
\hline
59&18.1040&17.75(Elkies,section-MW)&[59,27,16]\\
\hline
60&19.1039&19.04(Elkies,section-MW)&[60,28,16]\\
\hline
96&52.9003&52.078(\eta({\bf P}_{48q}))&[96,28,32]\\
\hline
136&100.157&100(\eta({\bf E}_8))&[136,57,32]\\
\hline
160&131.8847&127.4051(Analogous-Craig)&[160,27,64]\\
\hline
\end{array}
\]
\end{table}

\newpage

\section{Possible new lattices denser than Mordell-Weil lattices}

We prove the following result.

{\bf Proposition 5.1.} {\em Let $p$ be a prime satisfying $p \equiv
5$ $mod$ $6$. If there was a binary linear $[2p-2,
\frac{7p-5}{6}-\lceil\frac{p-11}{12}\cdot log_2 p \rceil, \geq
\frac{2(p+1)}{3}]$ code, then we could construct  a new lattice of dimension $
2p-2$ with center density larger than the center density
$\frac{((p+1)/12)^{p-1}}{p^{(p-5)/6}}$ of Mordell-Weil lattice in
this dimension.}
\begin{proof}
Set $m=\lceil\frac{p+1}{12}\rceil$. From Bertrand 's postulate there exists
a prime $l$ between $2p-1$ and $4p-2$. Applying Theorem 2.4
to analogous Craig lattice ${\bf A}_{2p-2}^{(m,l)}$ and the code in
the condition, we get the lattice.
\end{proof}

The center density of the dimension $52$ Mordell-Weil lattice is
$\frac{(5/2)^{26}}{2\cdot 53^4} \approx 2^{10.4578}$ from Theorem
1.2 \cite{Shi}(p.933 of \cite{Shi}). Craig lattice ${\bf A}_{52}^{(6)}$ has
its center density $\frac{3^{26}}{53^{5.5}}\approx 2^{9.7045}$.
Applying  Theorem 2.4 to this Craig lattice and the binary linear
$[52, 1, 52]$ code we get a lattice with center density
$2^{10.7045}$. This is a lattice with center
density slightly larger than the center density of Mordell-Weil
lattice from Theorem 1.2 of \cite{Shi}. It should be indicated that in the refinement of \cite{Flo} the least dimension is $57$. The dimension $140$ Mordell-Weil lattice has center density
$\frac{6^{70}}{71^{11}} \approx 2^{113.31}$.  Applying Theorem 2.4
to the binary $[140,50,32]$  code and the analogous Craig lattice
${\bf A}_{140}^{(4,151)}$ we get a lattice with center density
$2^{94.6656}$, which is slightly less than the center density $2^{97}$ of the non-lattice packing constructed
 in \cite{BiEd}. If there was a binary linear $[140,69,32]$ code (see
\cite{Gra}), there would be a possible new lattice of rank 140 with
center density at least $2^{114.6656}$, which is denser than Mordell-Weil lattice in this dimension .\\

{\bf Proposition 5.2.} {\em If there was binary linear $[128, 59,
32]$ code, then a new lattice of rank $128$ with center density
larger than the center density $2^{97.40}$ of Mordell-Wel lattice
could be constructed from Theorem 2.4.}
\begin{proof}
Applying Theorem 2.4 to ${\bf A}_{128}^{(4,131)}$ and the possible
binary linear $[128, 59,32]$ code or to ${\bf A}_{128}^{(6, 131)}$,
we could get the new lattice.
\end{proof}

From table in \cite{Gra} there exists a binary linear $[128,43,32]$
code we have a dimension $128$ lattice with center density
$2^{83.1784}$. The present upper bound for the minimum distance of
linear binary $[128, 59]$ code is $32$, but people do not know
whether this code exists or not. \\

{\bf Proposition 5.3.} {\em If there was one of the binary linear
$[256, 99, 64]$ code, $[256,74,80]$ code, $[256, 136, 48]$ code and
$[256,56,96]$ code, then a new 256 dimension lattice with center
density larger than the center density $2^{294.8}$ of Mordell-Weil
lattice could be constructed from Theorem 2.4.}
\begin{proof}
Applying Theorem 2.4 to analogous Craig lattice ${\bf A}_{256}^{(m,
257)}$, where $m=8,10,6,12$, and the possible binary linear code we
could get the lattice.
\end{proof}

From table in \cite{Gra} the present upper bound for the minimum
Hamming weight of a binary linear $[256,99]$ code is 74 and a binary
linear $[256,99, 48]$ code exists, this lead to a lattice with
center density $2^{257.8492}$. We do not know whether a linear
binary $[256,99,64]$ code exists or not. For the remaining cases, we
can find binary linear $[256, 74, d_{74}=56]$ code exists and the
present upper bound for $d_{74}$ is $85$, people do not know whether
a binary linear $[256,74,80]$ code exists or not. Similarly $32 \leq
d_{136} \leq 54$, people do not know whether a binary linear $[256,
136, 48]$ code exists or not. $68 \leq d_{56}
\leq 96$, people do not know whether a binary linear $[256, 56, 96]$ code exists or not.\\

For the dimension 508 case we have the following result.\\

{\bf Proposition 5.4.} {\em  If there was one of the linear
$[254,26,112]$ code over ${\bf F}_2$, $[254,36,104]$code over
 ${\bf F}_2$, $[254, 47, 96]$ code over ${\bf F}_2$, $[254,78,80]$ code over ${\bf F_2}$,
 $[169,18,104]$
 code over ${\bf F}_4$, $[169,24,96]$ code over ${\bf
F}_4$, $[169,39,80]$ code over ${\bf F}_4$, $[127,16,96]$ code over
${\bf F}_8$, $[127, 26, 80]$ code over ${\bf F}_8$, $[127, 42, 64]$
code over ${\bf F}_8$, then we could construct a new lattice in dimension
$508$ with center density larger than the center density
$2^{745.62}$ of Mordell-Weil lattice in this dimension(see
\cite{Shi} p.934).}
\begin{proof}
If there was a linear $[169, 24, 96]$ code over ${\bf F}_4$,
a linear binary $[507,48,\\192]$ code could be constructed as the
concatenated code with the binary linear $[3,2,2]$ code. From
Theorem 2.4 we could have a new lattice with center density $2^{48}
\delta({\bf A}_{508}^{(24, 509)}) \approx
2^{48+699.2897}=2^{747.2897}$. The other cases can be checked
similarly.
\end{proof}

From the table \cite{Gra} a linear $[169, 24, 85]$ code over ${\bf
F}_4$ has been constructed and the present upper bound for the
linear $[169, 24, d_{24}]$ code over ${\bf F}_4$ is $104$. People do
not know whether  a $[169, 24, 96]$
code over ${\bf F}_4$ exists or not.\\

If there was one of the binary linear $[512, 353, 64]$ code,
$[512,289, 80]$ code, $[512,441,48]$ code, $[512,239,96]$ code and
$[512,170,128]$ code, then a new lattice of rank $512$ with center
density larger than $2^{797.12}$ could be constructed from Theorem
2.4. We can apply Theorem 2.4 to the analogous Craig lattice ${\bf
A}_{1024}^{(m,1031)}$, where $m=10,12,16,24,30,32$, and the
following possible codes we could get the possible new lattice. If
there was one of the binary linear $[1024,925, 80]$ code,
$[1024,810,96]$ code, $[1024,598,128]$ code, $[1024,419,192]$ code,
$[1024,314,240]$ code and $[1024,286,256]$ code, then a new dense
lattice of rank $1024$ with center density larger than $2^{2018.2}$
could be constructed from Theorem 2.4.  In the case of of dimension
$2048$, if there was one of the binary linear $[2048,1479,192]$
code, $[2048,1216,240]$ code, $[2048,1142,256]$ code,
$[2048,719,384]$ code, $[2048,522,480]$ and $[2048,471,512]$ code,
then a new dense lattice of rank $2048$ with center density larger
than $2^{4891}$ could be constructed from Theorem 2.4. In the case
of dimension 4096, applying Theorem 2.4 to the analogous Craig
lattice ${\bf A}_{4096}^{(m,4099)}$, if there was one of the binary
linear $[4096,2708,384]$ code, $[4096,2192,480]$ code,
$[4096,2050,512]$ code, $[4096,865,960]$ and $[4096,770, 1024]$
code, then a new lattice of rank $4096$ with center density larger
than $2^{11527}$ could be constructed from Theorem 2.4.\\

For many lengthes in the range $128-256$ our knowledge about the
linear codes over ${\bf F}_2,{\bf F}_4,{\bf F}_8$ is not sufficient
to determine whether these {\em nice codes} in the table in
\cite{Gra} can be constructed or not. The motivation for the past
works on long binary codes are mainly from the construction of
efficient McEliece public key cryptosystem(for example, see
\cite{Can}). This is the first time  establishing an intimate
relation between long linear binary codes and the known densest
Mordell-Weil lattices in high dimensions. In the following table 4 some possible new lattices denser than
Mordell-Weil lattices are listed.

\begin{table}[ht]
\caption{}\label{eqtable}
\renewcommand\arraystretch{1.5}
\noindent\[
\begin{array}{|c|c|c|c|}
\hline
dimension&possible-log_2 \delta&known-densest&condition\\
\hline
128&98.3831&97.40(MW)&[128,59,32]\\
\hline
140&114.6656&113.31(MW)&[140,69,32]\\
\hline
164&148.1570&147.3318(MW)&[164,92,32]\\
\hline
164&147.3596&147.3318(MW)&[164,58,48]\\
\hline
176&165.8067&165.1474(MW)&[176,104,32]\\
\hline
176&166.3191&165.1474(MW)&[176,68,48]\\
\hline
200&194.9761&194.2188(MW)&[200,122,32]\\
\hline
200&195.0917&194.2188(MW)&[200,53,64]\\
\hline
212&221.4932&221.4145(MW)&[212,68,64]\\
\hline
224&241.3005&241.0012(MW)&[224,76,64]\\
\hline
256&294.958&294.8(MW)&[256,99,64]\\
\hline
256&295.15&294.8(MW)&[256,74,80]\\
\hline
256&294.8492&294.8(MW)&[256,136,48]\\
\hline
256&294.8156&294.8(MW)&[256,56,96]\\
\hline
272&323.1472&323.0536(MW)&[272,112,64]\\
\hline
380&525.4662&525.1006(MW)&[380,133,96]\\
\hline
452&671.0404&670.4412(MW)&[452,130,128]\\
\hline
508&747.2897&745.62(MW)&[508,48,192]\\
\hline
512&797.3117&797.12(MW)&[512,353,64]\\
\hline
692&1200.4738&1199.8554(MW)&[692,309,128]\\
\hline
716&1260.9065&1260.7960(MW)&[716,331,128]\\
\hline
1024&2018.2944&2018.2(MW)&[1024,286,256]\\
\hline
1436&3112.5083&3111.8561&[1436,571,256]\\
\hline
2048&4891.9666&4891(MW)&[2048,471,512]\\
\hline
4096&11527.8215&11527(MW)&[4096,770,1024]\\
\hline
\end{array}
\]
\end{table}

\newpage

\section{Lattice sphere packings denser than Mordell-Weil lattices}

{\bf Theorem 6.1.} {\em There exists a binary linear $[4096, 772, 1024]$ code. We can construct a lattice sphere
packing in dimension $4096$ with center density at least $2^{11529}$, which is denser than the Mordell-Weil
lattice in this dimension.}\\

\begin{proof}
From Gilbert-Varshamov bound if $V(4096,1023)=\Sigma_{i=0}^{1023} \displaystyle{4096 \choose i} < 2^{4097-k}$,
the linear binary $[4096,k,1024]$ code exists(\cite{Lint}). From the
inequality $\Sigma_{i=0}^r \displaystyle{n \choose i} <2^{nH(r/n)}<2^{H(\frac{1}{4})n}$, where $H(x)$ is the binary entropy function(\cite{Lint} p.21), then $V(4096,1023)<3324$. We get the conclusion.
\end{proof}

{\bf Lemma 6.2.} {\em The linear binary code of length $8n$, dimension $[(6log_2 3-8)n]$ and minimum distance $2n$ exists.}\\

\begin{proof}
Let $V(n,r)=\Sigma_{i=0}^r \displaystyle{n \choose i}$. From the Gilbert-Varshamov bound(\cite{Lint}), if $V(n,d-1) < 2^{n-k+1}$, then a linear binary code with parameter $[n,k,d]$ exists. From Theorem 1.4.5 of \cite{Lint} page 21 $V(8n, 2n-1) < 2^{8H(\frac{1}{4})n}$. The conclusion follows directly.
\end{proof}

Actually we have the following result.\\

{\bf Theorem 6.3.} {\em Let $p$ be a prime number satisfying $p \equiv 5$ $mod$ $6$ and $1667 \leq p \leq 2039$. We can construct
new lattice sphere packing in dimension $n=2p-2$ with center density better than the $\frac{((p+1)/12)^{p-1}}{p^{(p-5)/6}}$ of the corresponding Mordell-Weil lattice in dimension $n$.}

\begin{proof} The center density of Mordell-Weil lattice in dimension $2p-2$ is $\frac{((p+1)/12)^{p-1}}{p^{(p-5)/6}}$. The analogous Craig lattice ${\bf A}_{2p-2}^{([\frac{p-1}{16}],2p+t)}$, where $t$ is a positive integer such that $2p+t$ is a prime number, has its center density $\frac{[\frac{p-1}{32}]^{p-1}}{(2p+t)^{[\frac{p-1}{16}]-\frac{1}{2}}}$. It can be checked that $(2p+t)<2^{1.001}p$. On the other hand there exists a $[2p-2, [0.3776(p-1)],\frac{p-1}{2}]$ linear binary code, the conclusion follows directly.
\end{proof}

\remark The construction of the wanted linear binary codes in Theorem 6.3 is proved from the Gilbert-Varshamov bound. From the point view of coding theory, it seems quite possible that better linear binary codes exist. If that was true, the method in this paper would give lattice sphere packings in more dimensions better than Mordell-Weil lattices. It would be desirable that the codes in Theorem 6.1 and Theorem 6.3 could be constructed explicitly.
\endremark

In the following table 5 we list some new lattices denser than Mordell-Weil lattices.\\

\begin{table}[ht]
\caption{}\label{eqtable}
\renewcommand\arraystretch{1.5}
\noindent\[
\begin{array}{|c|c|c|}
\hline
dimension&new-log_2 \delta&known-densest\\
\hline
3332&8913&8897.0184(Moedell-Weil)\\
\hline
3956&11035&10969.9654(Mordell-Weil)\\
\hline
3992&11159&11099.6432(Mordell-Weil)\\
\hline
4004&11208&11130.5560(Mordell-Weil)\\
\hline
4052&11370&11294.2234(Mordell-Weil)\\
\hline
4076&11455&11375.6625(Mordell-Weil)\\
\hline
4096&11529&11527(Mordell-Weil)\\
\hline
\end{array}
\]
\end{table}

\newpage

\section{New denser lattices of high dimensions in the range 4098-8640}

The known densest sphere packing in dimension $4098$ is the Craig
lattice ${\bf A}_{4098}^{(246,4099)}$ with center density
$2^{11279}$(\cite{CS1}, page 17, Table 1.3). By using
the Craig lattice ${\bf A}_{4098}^{(128,4099)}$ ($4099$ is a prime,
\cite{prime}) and the linear binary $[4098, 773, 1024]$ code from
Lemma 6.2, we can construct a new dense lattice of rank $4098$ with
center density at least $2^{11536}$. In dimension $4104$ the lattice
$\eta({\bf \Lambda}_{24})$ has the center density
$2^{11400}$(\cite{CS1}, page 242, Table 8.7). By using the analogous Craig
lattice ${\bf A}_{4104}^{(128, 4111)}$($4111$ is a prime number) with
center density at least $2^{10780}$ and the
linear binary $[4104, 774, 1026]$ code, we can construct a new dense
lattice of rank $4104$  with center density at least $2^{11554}$. In
dimension $4124=2 \cdot 2063-2$, where $p=2063$ is a prime
satisfying $2063 \equiv 5$ $mod$ $6$, the Mordell-Weil lattice
sphere packing has the center density
$\frac{172^{2062}}{2063^{343}}\approx 2^{11537.1837}$. The analogous
Craig lattice ${\bf A}_{4124}^{(128,4127)}$($4127$ is a prime
number) has center density at least $2^{10840}$. From Lemma 6.2
there exists a linear binary $[4124, 778, 1031]$ code, we can construct a new denser lattice sphere packing with center density $2^{11618}$.\\

We have the following result.\\

{\bf Theorem 7.1.} {\em For each dimension $N=24n$ in the range $4104-8640$, we can
construct a new lattice sphere packing which is denser than the rank
$N$ child lattice of the Leech lattice $\eta({\bf \Lambda}_{24})$.}\\

\begin{proof}
We take the analogous Craig lattice ${\bf A}_{24n}^{([\frac{3n}{4}],l)}$, where $l$ is smallest prime number bigger than or equal to $24n+1$. From Lemma 6.2, there exists a linear binary $[24n, [4.5312n], 6n]$ code. Then we can construct a lattice with center density $2^{[4.5312n]} \cdot \frac{([\frac{3n}{8}])^{12n}}{l^{[3n/4]-\frac{1}{2}}}$. The conclusion follows from the direct computation \end{proof}

\remark It would be desirable if the codes(the existence is proved by Lemma 6.2 ) used in Theorem 7.1 could be constructed explicitly.
\endremark

In dimension $16392$, the known densest sphere packing is the child of the Leech lattice by Bos-Conway-Sloane construction with center density $2^{61608}$(\cite{CS1}, page 17). We can have a lattice sphere packing in dimension $16392$ with center density $2^{61497}$ from Theorem 2.4 and Lemma 6.2, which is less than the center density of the above known lattice sphere packing in dimension $16392$.\\

Here is the table 6 of new lattices of dimensions in the range $4098-16380$, which are better than the records in \cite{CS1}, page 17, Table 1.3.\\

\begin{table}[ht]
\caption{}\label{eqtable}
\renewcommand\arraystretch{1.5}
\noindent\[
\begin{array}{|c|c|c|}
\hline
dimension&new-log_2 \delta&old-record\\
\hline
4098&11536&11279(Craig)\\
\hline
4104&11554&11400(\eta({\bf \Lambda}_{24}))\\
\hline
4124&11618&11537.1837(Mordell-Weil)\\
\hline
8184&26823&26712(\eta({\bf \Lambda}_{24}))\\
\hline
8190&26915&26154(Craig)\\
\hline
8208&26953&26808(\eta({\bf \Lambda}_{24}))\\
\hline
16380&61419&59617(Craig)\\
\hline
\end{array}
\]
\end{table}

\section{New dense lattice sphere packings of odd dimensions}

In many previous constructions of dense lattice sphere packing in
the Euclid space of dimension $n$, $n$ is always assumed to be an
even number. Craig lattice is of dimension $p-1$, where
$p$ is a prime number. The Mordell-Weil lattices and the children of
the Leech lattice in the  Bos-Conway-Sloane construction are of even
dimensions. Our construction Theorem 2.2 and 2.4 can be used for
each dimension and the center density of the constructed lattice
sphere packing depends "continuously" on the dimension. Thus we can
have quite good lattice sphere packings in the "missing" or odd
dimensions. In the following table 7, new dense lattice sphere
packings in moderate dimensions $149-159$ and $183-193$ are listed.
It should be noted that in these dimensions, no previous
constructions like Craig, Bos-Conway-Sloane and Mordell-Weil can be
applied, except in dimension $150,156,190$ and $192$,
there are known Craig lattices and their refinements(\cite{Flo}, \cite{CS1}, page 17). In the dimensions $150, 156, 190$ and $192$, the refinements of Craig lattices in \cite{Flo} are slightly denser.\\

\begin{table}[ht]
\caption{}\label{eqtable}
\renewcommand\arraystretch{1.5}
\noindent\[
\begin{array}{|c|c|c|c|}
\hline
dimension&new-log_2 \delta&dimension&new-log_2\\
\hline
149&112.3048&183&158.4505\\
\hline
150&114.06(<[18])&184&160.0355\\
\hline
151&113.7424&185&161.6205\\
\hline
152&115.2811&186&163.2055\\
\hline
153&116.8248&187&164.7905\\
\hline
154&118.3685&188&166.3755\\
\hline
155&119.9122&189&167.9605\\
\hline
156&121.4559(<[18])&190&169.5455(<[18])\\
\hline
157&122.1067&191&171.1305\\
\hline
158&123.6504&192&172.44(<[18])\\
\hline
159&125.1941&193&173.5188\\
\hline
\end{array}
\]
\end{table}

\remark We observe that the above lattice sphere packings are much
better than the lattice sphere packings from Minkowski-Hlawka
Theorem. New dense lattice sphere packings can also be constructed
from Theorem 2.4 for dimensions $n=241,...,251$.
\endremark

In table 8 we list some new lattice sphere packings in the dimension
$p,p-2$ and $p+2$ where $p$ is a prime number. The center density is
quite close to the center density of  the Craig lattice  in the
dimension $p-1$(\cite{CS1}, page 17, Table 1.3). In table 9 some new
lattice sphere packings in the dimensions $2p-3$, where $p$ is a
prime number satisfying $p \equiv 5$ $mod$ $6$, are listed. We also
list some new lattice sphere packings in dimensions
$24n'-1$ in table 10.\\

\begin{table}[ht]
\caption{}\label{eqtable}
\renewcommand\arraystretch{1.5}
\noindent\[
\begin{array}{|c|c|c|c|c|c|}
\hline
dimension&log_2 \delta&dimension&log_2 \delta&dimension&log_2 \delta\\
\hline
87&40.4835&89&42.5005&91&43.9503\\
\hline
149&112.3048&151&115.0103&153&117.4377\\
\hline
179&153.5829&181&155.3909&183&158.4792\\
\hline
189&167.7417&191&170.5800&193&173.1791\\
\hline
507&741.1263&509&744.4672&511&748.8247\\
\hline
\end{array}
\]
\end{table}

\begin{table}[ht]
\caption{}\label{eqtable}
\renewcommand\arraystretch{1.5}
\noindent\[
\begin{array}{|c|c|c|}
\hline
dimension&new-log_2 \delta&known-densest(dim=n+1)\\
\hline
3331&8910.1498&8897.0184(Moedell-Weil)\\
\hline
3955&11032.7450&10969.9654(Mordell-Weil)\\
\hline
3991&11156.0229&11099.6432(Mordell-Weil)\\
\hline
4003&11212.0171&11130.5560(Mordell-Weil)\\
\hline
4051&11365&11294.2234(Mordell-Weil)\\
\hline
4075&11452&11375.6625(Mordell-Weil)\\
\hline
4095&11526&11527(Mordell-Weil)\\
\hline
\end{array}
\]
\end{table}

\begin{table}[ht]
\caption{}\label{eqtable}
\renewcommand\arraystretch{1.5}
\noindent\[
\begin{array}{|c|c|c|}
\hline
dimension&new-log_2 \delta&old-record(dim=n+1)\\
\hline
4097&11533&11279(Craig)\\
\hline
4103&11551&11400(\eta({\bf \Lambda}_{24}))\\
\hline
4123&11615&11537.1837(Mordell-Weil)\\
\hline
8183&26819&26712(\eta({\bf \Lambda}_{24}))\\
\hline
8189&26911&26154(Craig)\\
\hline
8207&26949&26808(\eta({\bf \Lambda}_{24}))\\
\hline
16379&61415&59617(Craig)\\
\hline
\end{array}
\]
\end{table}

\end{document}